\documentclass[12pt]{amsart}

\usepackage{amsfonts}
\usepackage{amsthm}
\usepackage{amsmath}
\usepackage{amssymb}

\newcommand{\reali}{\mathbb{R}}
\newcommand{\complessi}{\mathbb{C}}
\newcommand{\nat}{\mathbb{N}}
\newcommand{\ident}{{\mathchoice {\rm 1\mskip-4mu l} {\rm 1\mskip-4mu l}
{\rm 1\mskip-4.5mu l} {\rm 1\mskip-5mu l}}}
\newcommand{\pr}{\mathbb{P}}

\newtheorem{teo}{Theorem}[section]
\newtheorem{lem}[teo]{Lemma}
\newtheorem{cor}[teo]{Corollary}
\newtheorem{rem}[teo]{Remark}

\def\eqref#1{(\ref{#1})}

\begin{document}

\title{A note on compact Markov operators}

\author[Fabio Zucca]{Fabio Zucca \\
Dipartimento di Matematica \\
Politecnico di Milano\\
Piazza Leonardo da Vinci 32, 20133 Milano, Italy.\\
zucca@mate.polimi.it}

\renewcommand{\subjclassname}{\textup{2000} Mathematics Subject Classification}

\begin{abstract}
The analytic properties of the Markov operator associated to a
random walk are common tools in the study of the behaviour and
some probabilistic features related to the walk. In this paper
we consider a class of Markov operators which generalizes the 
class of  compact Markov operators and we study some probabilistic
properties of the associated random walk.
\end{abstract}

\maketitle



\section{Basic definitions}

\label{sec:basic}

Let $(X,P)$ be an irreducible, random walk on the
state space $X$ which is at most countable. We suppose that
the (usually infinite) stochastic matrix $P$
describes a Markov chain
$\{Z_n\}_{n\in \nat}$ defined
on a probability space $(\Omega,\mathcal{F}, \mathbb P)$ with
transition probabilities
$p(x,y):={\mathbb P}[Z_{n+1}=y|Z_n=x]$
homogeneous in time.
Besides we consider the $n$-step transition probabilities
$\{p^{(n)}(x,y)\}_{x,y\in X}$ which represent the stochastic
matrix associated to the $n$-th convolution power of $P$.

The {\it Markov operator} associated to the random walk is defined
as follows
\begin{equation}\label{operatorep}
\begin{split}
D(P)&:=\left\{f:X \rightarrow \reali :
\sum_{y \in X} p(x,y) |f(y)| <+\infty, \ \forall x\in X \ \right\}, \\
(Pf)(x)&:=\sum_{y \in X} p(x,y) f(y), \qquad \forall f \in D(P), \ \forall
x \in X;\\\
\end{split}
\end{equation}
note that $\mathcal{D}(P) \supseteq l^\infty(X)$ and that
$P|_{l^\infty(X)}$
is a bounded linear operator from $l^\infty(X)$ into itself.

To explore the behaviour of the random walk $(X,P)$ and its main properties
we introduce the two generating functions
\[
\begin{split}
G(x,y|z)&=\sum_{n=0}^\infty p^{(n)}(x,y) z^n, \qquad
F(x,y|z)=\sum_{n=0}^\infty f^{(n)}(x,y) z^n \\
\end{split}
\]
where $\{f^{(n)}(x,y)\}_{x,y\in X}$ are the first time return probabilities, namely
\[
f^{(n)}(x,y)=\pr(Z_n=y, Z_i\not = y , \forall i=1,\ldots n-1|Z_0=x), \quad f^{(0)}(x,y)=0.
\]
Both the generating functions must be considered inside their circle of convergence in
$\complessi$.

An irreducible random walk $(X,P)$ is called {\it transient} if and only if there exists ($\Leftrightarrow$ for any) $x\in X$
such that $F(x,x)<1$ and {\it recurrent} otherwise. Among the recurrent random walks
we distinguish the class of {\it positive recurrent} and {\it null recurrent} depending
on whether $\overline{\tau}_x:=
\sum_{n=1}^\infty n f^{(n)}(x,x) <+\infty$ for some ($\Leftrightarrow$ for any) $x\in X$
or not.

We note that positive recurrence is a strong assumption: for instance
if $(X,P)$ is the simple random walk on a infinite, locally finite, non-oriented,
connected graph $(X,E)$, then it is not positive recurrent.
Indeed it is easily reversible
with reversibility measure given by
$
m(x):= \# \{y: (x,y) \in E\}
$, which is clearly infinite.
According to Theorem 1.18 of \cite{Woess2},
if an irreducible Markov chain is recurrent, then it admits a unique (up to multiplication)
stationary measure and this one is finite if and only if the walk is positive
recurrent. Since a reversibility measure is stationary,
if the walk were positive recurrent, then $m$ should be finite.

The importance of this class of random walks is highlighted by Theorem 1.18
of \cite{Woess2} (see also Theorem 3.2 of \cite{Zucca4}).

\begin{rem}\label{positive}
We note
that $(X,P)$ is positive recurrent if and only if there exists
($\Leftrightarrow$ for all) $x\in X$,
\[
\lim_{z \rightarrow 1^-} \frac{1-F(x,x|z)}{1-z} <+\infty.
\]
Just take in mind
that the limit always exists
(finite or infinite)
due to the decomposition
\[
\lim_{z \rightarrow 1^-}\frac{1-F(x,x|z)}{1-z} =
\lim_{z \rightarrow 1^-}\frac{1-F(x,x|1)}{1-z}+
\lim_{z \rightarrow 1^-}\frac{F(x,x|1)-F(x,x|z)}{1-z},
\]
and, using well-known arguments, we have
\[
\begin{split}
\sum_{n=0} n f^{(n)}(x,x)&= \lim_{z \rightarrow 1^-}
\sum_{n=0} n f^{(n)}(x,x) z^n \\
& =\lim_{z \rightarrow 1^-} F^\prime(x,x|z) =
\lim_{z \rightarrow 1^-} \frac{F(x,x|1)-F(x,x|z)}{1-z};
\end{split}
\]
(note that the last equality holds also if the limit
is $+\infty$).
\end{rem}

\section{Compact Markov operators}

\label{sec:compact}

In this section we want to study the behaviour of a random walk whose associated 
Markov operator satisfying equation~(\ref{eq:compact2}) below. In particular
we study compact Markov operators.


We recall here the characterization of a compact operator
(with non-negative
matrix elements)
defined by equation \eqref{operatorep}
(see \cite{Zucca4}, Theorem 2.2). 

\begin{teo}\label{compact1}
Let $X$ be a countable set and
let $P$ be a transition operator
on $X$ with non negative elements, satisfying the condition
$\sup_{x \in X} \sum_{y \in X} p(x,y) <+\infty$. Then $P$ is a bounded,
linear operator from $l^\infty(X)$ into itself; moreover $P:l^\infty \mapsto l^\infty$
is compact
if and only if for any given $\epsilon>0$ there exists a finite subset
$A_\epsilon \subset X$ such that
$\sup_{x \in X} \sum_{y \in X\setminus A_\epsilon} p(x,y) < \epsilon$.
\end{teo}

The next Theorem is the main result of this section: its corollary
enhances Proposition 2.3 of \cite{Zucca4}.

\begin{teo}\label{compact2}
Let $(X,P)$ be an irreducible Markov chain, and suppose that there exists
$\epsilon \in (0,1)$ and a finite subset $A\subset X$ such that
\begin{equation}\label{eq:compact2}
\sup_{x \in X} \sum_{y \in X\setminus A} p(x,y) < \epsilon;
\end{equation}
then $(X,P)$ is positive recurrent.
\end{teo}

Before looking at the proof, we want to understand what equation \eqref{eq:compact2} implies from the
point of view of the walker.

Let us consider the preadjoint map $P_*:l^1(X) \rightarrow l^1(X)$ acting as
\[
P_*\nu(y)\equiv \nu P(y) := \sum_{x \in X} \nu(x)p(x,y),\qquad \forall y\in X.
\]
This is a mass-preserving map, indeed $\sum_{y \in X} \nu P(y) = \sum_{x\in X} \nu(x)$;
moreover $\nu \geq 0$ implies $P_* \nu \geq 0$.
If $\nu$ is the probability distribution of the
position of the walker at a certain time, then
$P_* \nu$
is
the probability distribution after one step.

In term of this evolution map, 
equation~\eqref{eq:compact2} is equivalent to the existence of 
a finite subset $A$ such that, given any
probability distribution $\nu$, the probability distribution after one step satisfies
$P_* \nu(A) \geq 1-\epsilon$
(or, equivalently, $P_*^n \nu(A)\geq 1-\epsilon$ for any $n \in \nat^*$).

Since from the Law of large numbers, for any given $x \in X$, $\pr$-a.c.
\[
\lim_{n \rightarrow +\infty} \frac{1}{n} \sum_{i=1}^n \ident_{\{x\}}(Z_n)=
\begin{cases}
0 & \text{ in the transiente or null-recurrent case}\\
\frac{1}{\overline\tau}_x & \text{ in the positive recurrent case}
\end{cases}
\] 
and since in the positive recurrent case, $1/\overline{\tau}_x$ represents the unique
stationary probability measure
(see \cite{Chung1}, Section I.7, Theorem~1), hence the walker will pass (asymptotically) at least $1-\epsilon$ of
its time in $A$.

\begin{proof} {\it (of Theorem \ref{compact2}).}
Let $A$ and $\epsilon$ satisfying equation \eqref{eq:compact2}.

Let us note that for any given $x \in X$ we have
\[
\sup_{x \in X} \sum_{y \in X\setminus A} p^{(n)}(x,y) < \epsilon,
\]
indeed for any given $n \in \nat^*$
\[
\pr(Z_n \in A |Z_0=x)= \sum_{z \in X}
\pr(Z_n \in A |Z_{n-1}=z)\pr(Z_{n-1}=z |Z_0=x) \ge 1-\epsilon.
\]
Let $x_0 \in X \setminus A$ be fixed and rewrite the previous equation
as
\[
\sum_{x\in A} p^{(n)}(x_0,x) \ge 1-\epsilon, \qquad \forall n \in \nat^*.
\]
This implies
\[
\sum_{x\in A} G(x_0,x|z) \ge (1-\epsilon) \left ( \sum_{j=1}^\infty z^j \right )=
\frac{z(1-\epsilon)}{1-z}.
\]
Since $G(x,y|z)=\delta_{xy}+F(x,y|z)G(y,y|z) =\delta_{xy}+F(x,y|z)/(1-F(y,y|z))$,
the previous inequality becomes
\[
\sum_{x\in A} \frac{F(x_0,x|z)}{1-F(x,x|z)} \ge \frac{z(1-\epsilon)}{1-z}.
\]
Now, taking in mind the usual position $1/\infty=0$, since $A$ is finite,
\[
\begin{split}
0<
1-\epsilon \le \lim_{z \rightarrow 1^-} \sum_{x\in A}
\frac{(1-z) F(x_0,x|z)}{1-F(x,x|z)} =
\sum_{x \in A} \frac{F(x_0,x|1)}{\lim_{z \rightarrow 1^-}(1-F(x,x|z))/(1-z)}.
\end{split}
\]
This easily implies the existence of $x_1 \in A$ such that
\[
\lim_{z \rightarrow 1^-}\frac{1-F(x_1,x_1|z)}{1-z} < +\infty \\
\]
hence, by Remark~\ref{positive}, $(X,P)$ is positive recurrent.
\end{proof}

\begin{cor}\label{compact2cor}
Let $(X,P)$ be an irreducible Markov chain, such that $P$ is compact;
then $(X,P)$ is positive recurrent.
\end{cor}

We emphasize that the Markov operator associated to a
positive recurrent random walk needs not to be compact as the following
example shows.
Take $X=\nat$ and define the transition probabilities as follows:
\[
\begin{cases}
\begin{matrix}
p(0,1)=1, \hfill& \\
p(n,n+1)=1-p,  \hfill & \forall n \in \nat^*,  \hfill \\
p(n+1,n)=p,  \hfill & \forall n \in \nat,  \hfill \\
0  \hfill& \text{otherwise,}  \hfill \\
\end{matrix}
\end{cases}
\]
where $p\in (0,1)$.
The first time return probabilities generator function $F$ can be easily calculated for  $x=y=0$
as
\[
F(0,0|z)=\frac{2pz^2}{1+\sqrt{1-4z^2p(1-p)}};
\]
the corresponding random walk is transient if $p\in (0,1/2)$, null recurrent if $p=1/2$ and
positive recurrent if $p\in (1/2,1)$, but equation~\eqref{eq:compact2} does not hold,
hence the Markov operator is always non compact.

\section{Some estimates}

\label{sec:estimate}

Let us define the time of the first return ont the vertex $x\in X$
as ${\rm T}_x:=\inf\{n \ge 1: Z_n=x\}$; 
in the recurrent case we have
\[
\overline{\tau}_x=
{\mathbb E}[{\rm T}_x|Z_0=x]= \lim_{z \rightarrow 1^-} \frac{1-F(x,x|z)}{1-z}
\equiv \lim_{z \rightarrow 1^-} (1-z)G(x,x|z).
\]
Moreover if equation \eqref{eq:compact2} holds, we have that $F(x,y)=1$ for any
$x,y \in X$ and
\[
1 \ge
\sum_{x\in A} \frac{1}{\overline{\tau}_x} \ge 1-\epsilon
\]
which implies
\[
\min_{x \in A}
\overline{\tau}_x \le \frac{{\rm card}(A)}{1-\epsilon}.
\]

Besides for the first time entrance in $A$, ${\rm T}_A:=\inf \{ n>0: Z_n \in A\}$,
the following hold
\[
{\mathbb P}({\rm T}_A\ge n| Z_0=x) \le \epsilon^{n-1}, \qquad
{\mathbb E}[{\rm T}_A|Z_0=x] \le \frac1{1-\epsilon}.
\]

In the 
reversible case, it is possible to find lower bounds for
the $n$-step transition
probabilities $p^{(n)}(x,x)$ and their generating function.
In this case the reversibility measure $m$ satisfies $m(x)\propto 1/{\overline \tau}_x$.

\begin{lem}\label{unilow}
Let $(X,P)$ be a reversible random walk, $m$ a reversibility measure
and $x\in \Gamma$ a fixed vertex.
If there exists $n\in \nat^*$  and $A \subseteq X$ such that
\[
\sum_{y \in X \setminus A} p^{(n)}(x,y) \le \epsilon,
\]
then
\[
p^{(2n)}(x,x) \geq (1-\epsilon)^2 \frac{m(x)}{m(A)}.
\]
\end{lem}

The easy proof of this lemma is straightforward and we omit it.
By using this result one see immediately that, for any $x \in A$
\[
\begin{split}
G(x,x|z) & \ge \frac{(1-\epsilon)^2}{1- z^2} \frac{m(x)}{m(A)}, \qquad z \in [0,1).
\end{split}
\]




\end{document}